\documentclass[a4paper,12pt,english]{smfart}

\usepackage{amssymb}
\usepackage{url}
\RequirePackage{mathrsfs}
\DeclareSymbolFont{rsfscript}{OMS}{rsfs}{m}{b}
\DeclareSymbolFontAlphabet{\mathrsfs}{rsfscript}
\usepackage[utopia,expert]{mathdesign}
\usepackage{lscape}
\usepackage{calligra}
\usepackage[T1]{fontenc}

\usepackage{frcursive}

\usepackage{palatino}
\usepackage{rotating}
\usepackage{graphicx}

\usepackage{enumerate}

\usepackage{color}
\definecolor{shadecolor}{gray}{0.80}

\def\bfit{\bfseries\itshape}

\def\proofname{D\'emonstration}

\input xypic
\xyoption{all}
\xyoption{arc}

\def\thetheo{\thesection.\arabic{theo}}

\renewcommand\theequation{\thetheo}
\def\equat{\refstepcounter{theo}\begin{equation}}
\def\endequat{\end{equation}}

\renewcommand\thetable{\thetheo}

\renewcommand\thesection{\arabic{section}}
\renewcommand\thesubsection{\thesection.\Alph{subsection}}

\def\contentsname{Table des mati\`eres}
\def\listfigurename{Liste des figures}
\def\listtablename{Liste des tables}
\def\appendixname{Appendice}

\def\refname{R\'ef\'erences}

  \def\lG{{\mathfrak l}}  
    
    \def\NM{{\mathbb{N}}}

    \def\RM{{\mathbb{R}}}
  \def\sG{{\mathfrak s}}

  \def\bb{{\mathbf b}}

    \def\FC{{\mathcal{F}}}
    \def\GC{{\mathcal{G}}}

  \def\tb{{\mathbf t}}

    \def\WC{{\mathcal{W}}}

  \def\zb{{\mathbf z}}

\def\Crm{{\mathrm{C}}}

\def\Mrm{{\mathrm{M}}}

\def\b{\beta}
\def\g{\gamma}

\def\e{\varepsilon}
\def\ph{\varphi}

\def\L{\Lambda}

\def\o{\omega}
\def\O{\Omega}

\def\lamb{{\boldsymbol{\lambda}}}       
       \def\Lamt{{\tilde{\Lambda}}}

\DeclareMathOperator{\Bip}{{\mathrm{Bip}}}

\DeclareMathOperator{\Irr}{{\mathrm{Irr}}}

\DeclareMathOperator{\Cons}{{\mathrm{Cons}}}
\DeclareMathOperator{\Fam}{{\mathrm{Fam}}}
\DeclareMathOperator{\Cell}{{\mathrm{Cell}}}

\def\to{\rightarrow}
\def\longto{\longrightarrow}

\def\fonctio#1#2#3#4{\begin{array}{ccc}
{#1} & \longto & {#2} \\
{#3} & \longmapsto & {#4} 
\end{array}}

\def\vide{\varnothing}
\def\DS{\displaystyle}
\def\SS{\scriptstyle}

\def\fin{~$\blacksquare$}
\def\finl{~$\blacksquare$}

\def\lexp#1#2{\kern\scriptspace\vphantom{#2}^{#1}\kern-\scriptspace#2}
\def\le{\hspace{0.1em}\mathop{\leqslant}\nolimits\hspace{0.1em}}
\def\ge{\hspace{0.1em}\mathop{\geqslant}\nolimits\hspace{0.1em}}

\mathchardef\inferieur="321E
\mathchardef\superieur="321F

\def\eqna{\begin{eqnarray*}}
\def\endeqna{\end{eqnarray*}}

\def\itemth#1{\item[${\mathrm{(#1)}}$]}

\catcode`\@=11
\long\def\@car#1#2\@nil{#1}
\long\def\@first#1#2{#1}
\long\def\@second#1#2{#2}
\long\def\ifempty#1{\expandafter\ifx\@car#1@\@nil @\@empty
  \expandafter\@first\else\expandafter\@second\fi}
\catcode`\@=12

\renewcommand{\tocsection}[3]{%
  {\bf \indentlabel{\ifempty{#2}{}{\ignorespaces#1 #2.\quad}}#3}}

\def\boitegrise#1#2{\begin{centerline}{\fcolorbox{black}{shadecolor}{~
    \begin{minipage}[t]{#2}{\vphantom{~}#1\vphantom{$A_{\DS{A_A}}$}}
            \end{minipage}~}}\end{centerline}\medskip}

\theoremstyle{remark}

\theoremstyle{plain}

\def\calo{{\Crm\Mrm}}

\def\xyinj{\ar@{^{(}->}}
\def\xysur{\ar@{->>}}

\bigskip

\def\kl{{\mathrm{Lus}}}

\makeatletter
\def\hlinewd#1{%
\noalign{\ifnum0=`}\fi\hrule \@height #1 %
\futurelet\reserved@a\@xhline}
\makeatother

\newlength\epaisLigne

\usepackage{multirow}

\def\ref{{\mathrm{ref}}}

\begin{document}

\baselineskip=16pt
%\large\baselineskip=20pt
%\Large\baselineskip=24pt

\title{Constructible characters and ${\boldsymbol{b}}$-invariant}

\author{{\sc C\'edric Bonnaf\'e}}
\address{
Institut de Math\'ematiques et de Mod\'elisation de Montpellier (CNRS: UMR 5149), 
Universit\'e Montpellier 2,
Case Courrier 051,
Place Eug\`ene Bataillon,
34095 MONTPELLIER Cedex,
FRANCE} 

\makeatletter
\email{cedric.bonnafe@univ-montp2.fr}
\makeatother

% \author{{\sc Rapha\"el Rouquier}}
% 
% \address{Mathematical Institute, 
% University of Oxford, 24-29 St Giles', Oxford, OX1 3LB, UK}
% \email{rouquier@maths.ox.ac.uk}

%\makeatother

\keywords{Coxeter groups, $b$-invariant, families, cellular characters}

\subjclass{According to the 2000 classification: Primary 20F55; Secondary 20C08}

\date{\today}

\thanks{The author is partly supported by the ANR (Project No ANR-12-JS01-0003-01 ACORT)}

\begin{abstract} 
If $W$ is a finite Coxeter group and $\ph$ is a weight function, Lusztig has defined 
{\it $\ph$-constructible characters} of $W$, as well as a partition of the set of irreducible 
characters of $W$ into {\it Lusztig $\ph$-families}. We prove that every Lusztig $\ph$-family 
contains a unique character with minimal $b$-invariant, and that every $\ph$-constructible 
character has a unique irreducible constituent with minimal $b$-invariant. This 
generalizes Lusztig's result about {\it special characters} to the case where $\ph$ 
is not constant. This is compatible with some conjectures of Rouquier and the author about 
{\it Calogero-Moser families} and {\it Calogero-Moser cellular characters}.
\end{abstract}

\maketitle

\pagestyle{myheadings}

\markboth{\sc C. Bonnaf\'e}{\sc Constructible characters and $b$-invariants}

% \tableofcontents
% 
% \vskip1cm

Let $(W,S)$ be a finite Coxeter system and let $\ph : S \to \RM_{>0}$ 
be a {\it weight function} that is, a map such that $\ph(s)=\ph(t)$ whenever 
$s$ and $t$ are conjugate in $W$. Associated with this datum, G. Lusztig has 
defined~\cite[\S{22}]{lusztig} a notion of {\it constructible characters} of $W$: 
it is conjectured that a character is constructible if and only if it is the character 
afforded by a Kazhdan-Lusztig left cell (defined using the weight function $\ph$). 
These constructible characters depend heavily on $\ph$ so we will call them 
the {\it $\ph$-constructible characters of $W$}: the set of 
$\ph$-constructible characters will be denoted by $\Cons_\ph^\kl(W)$. 
We will also define a graph $\GC_{W,\ph}^\kl$ as follows: the vertices of $\GC_{W,\ph}^\kl$ 
are the irreducible characters and two irreducible characters $\chi$ and $\chi'$ are joined 
in this graph if there exists a $\ph$-constructible character $\g$ of $W$ such that 
$\chi$ and $\chi'$ both occur as constituents of $\g$. The connected components 
of $\GC_{W,\ph}^\kl$ (viewed as subsets of $\Irr(W)$) will be called the 
{\it Lusztig $\ph$-families}: the set of Lusztig $\ph$-families will 
be denoted by $\Fam_\ph^\kl(W)$. If $\FC \in \Fam_\ph^\kl(W)$, we denote by 
$\Cons_\ph^\kl(\FC)$ the set of $\ph$-constructible characters of $W$ all of whose irreducible 
components belong to $\FC$. 

On the other hand, using the theory of rational Cherednik algebras at $t=0$ and the geometry 
of the Calogero-Moser space associated with $(W,\ph)$, R. Rouquier and the 
author (see~\cite{note} and~\cite{cmcells}) 
have defined a notion of {\it Calogero-Moser $\ph$-cells} of $W$, 
a notion of {\it Calogero-Moser $\ph$-cellular characters} of $W$ 
(whose set is denoted by $\Cell_\ph^\calo(W)$) and a notion 
of {\it Calogero-Moser $\ph$-families} (whose set is denoted by 
$\Fam_\ph^\calo(W)$).

\bigskip

\begin{quotation}
\noindent{\bf Conjecture (see~\cite{note},~\cite{cmcells} and~\cite{gordon martino}).} 
{\it With the above notation,
$$\Cons_\ph^\kl(W)=\Cell_\ph^\calo(W)\quad\text{\it and}\quad \Fam_\ph^\kl(W)=\Fam_\ph^\calo(W)$$
for every weight function $\ph : S \to \RM_{>0}$.} 
\end{quotation}

\bigskip

The statement about families in 
this conjecture holds for classical Weyl groups 
thanks to a case-by-case analysis relying on~\cite[\S{22}]{lusztig} 
(for the computation of Lusztig $\ph$-families),~\cite{gordon martino} 
(for the computation of Calogero-Moser $\ph$-families in type $A$ and $B$) 
and~\cite{bellamy} (for the computation of the Calogero-Moser $\ph$-families in type $D$). 
It also holds whenever $|S|=2$ 
(see~\cite[\S{17}~and~Lemma~22.2]{lusztig} and~\cite[\S{6.10}]{bellamy these}). 

The statement about constructible characters is much more difficult to establish, 
as the computation of Calogero-Moser $\ph$-cellular characters is at that time 
out of reach. It has been proved whenever the Caloger-Moser space associated 
with $(W,S,\ph)$ is smooth~\cite[Theorem~14.4.1]{cmcells} (this includes the cases where $(W,S)$ 
is of type $A$, or of type $B$ for a large family of weight functions: in all these cases, 
the $\ph$-constructible characters are the irreducible ones). It has also been checked by 
the author whenever $|S|=2$ or $(W,S)$ is of type $B_3$ (unpublished).

Our aim in this paper is to show that this conjecture is compatible 
with properties of the $b$-invariant (as defined below). 
With each irreducible character $\chi$ of $W$ 
is associated its {\it fake degree} $f_\chi(\tb)$, 
using the invariant theory of $W$ (see for instance~\cite[Definition~1.5.7]{cmcells}). 
Let us denote by $b_\chi$ the valuation of $f_\chi(\tb)$: 
$b_\chi$ is called the {\it $b$-invariant} of $\chi$. 
Let $r_{\! \chi}$ denote the coefficient of $\tb^{b_\chi}$ in $f_\chi(\tb)$. 
In other words, 
$$r_{\! \chi} \in \NM^*\quad\text{and}\quad f_\chi(\tb) \equiv r_{\!\!\chi}\tb^{b_\chi} \mod \tb^{b_\chi + 1}.$$
% In other words, 
% $$f_\chi \in \NM^*\quad\text{and}\quad f_\chi(\tb) \equiv f_\chi~\tb^{b_\chi} \mod t^{b_\chi + 1}.$$
For instance, $b_1=0$ and $b_\e$ is the number of reflections 
of $W$ (here, $\e : W \to \{1,-1\}$ denotes the sign character). Also, $b_\chi = 1$ 
if and only if $\chi$ is an irreducible constituent of the canonical reflection 
representation of $W$. The following result is proved in~\cite[Theorems~9.6.1~and~12.3.14]{cmcells}:

\bigskip

\noindent{\bf Theorem CM.} {\it 
Let $\ph : S \to \RM_{>0}$ be a weight function. Then:
\begin{itemize}
\itemth{a} If $\FC \in \Fam_\ph^\calo(W)$, then there exists a unique $\chi_\FC \in \FC$ 
with minimal $b$-invariant. Moreover, $r_{\chi_\FC}=1$. 

\itemth{b} If $\g \in \Cell_\ph^\calo(W)$, then 
there exists a unique irreducible constituent $\chi_\g$ of $\g$ with minimal $b$-invariant. 
Moreover, $r_{\chi_\g}=1$. 
\end{itemize}}

\bigskip

The next theorem is proved in~\cite[Theorem~5.25~and~its~proof]{lusztig orange} 
(see also~\cite{lusztig special} for the first occurence of the {\it special} 
representations):

\bigskip

\noindent{\bf Theorem (Lusztig).} 
{\it Assume that $\ph$ is {\bfit constant}. Then:
\begin{itemize}
\itemth{a} If $\FC \in \Fam_\ph^\kl(W)$, then there exists a unique $\chi_\FC \in \FC$ 
with minimal $b$-invariant ($\chi_\FC$ is called the {\bfit special} character of $\FC$). 
Moreover, $r_{\chi_\FC}=1$.

\itemth{b} If $\g \in \Cons_\ph^\kl(\FC)$, then $\chi_\FC$ is an irreducible constituent 
of $\g$ (and, of course, among the irreducible constituents of $\g$, $\chi_\FC$ is the 
unique one with minimal $b$-invariant). Moreover, $\langle \g,\chi_\FC \rangle =1$. 
\end{itemize}}

\bigskip

It turns out that, for general $\ph$, there might exist Lusztig $\ph$-families $\FC$ 
such that no element of $\FC$ occurs as an irreducible constituent of {\it all} the 
$\ph$-constructible characters in $\Cons_\ph^\kl(\FC)$ (this already occurs in type $B_3$, 
and the reader can also check this fact in type $F_4$, using the tables 
given by Geck~\cite[Table~2]{geck f4}). Nevertheless, 
we will prove in this paper the following result, which is compatible with the 
above conjecture and the above theorems:

\bigskip

\noindent{\bf Theorem L.} 
{\it Let $\ph : S \to \RM_{>0}$ be a weight function. Then:
\begin{itemize}
\itemth{a} If $\FC \in \Fam_\ph^\kl(W)$, then there exists a unique $\chi_\FC \in \FC$ 
with minimal $b$-invariant. Moreover, $r_{\chi_\FC}=1$. 

\itemth{b} If $\g \in \Cons_\ph^\kl(W)$, then 
there exists a unique irreducible constituent $\chi_\g$ of $\g$ with minimal $b$-invariant. 
Moreover, $r_{\chi_\g}=1$ and $\langle \g,\chi\rangle =1$. 
\end{itemize}}

\bigskip

The proof of Theorem~CM is general and conceptual, while 
our proof of Theorem~L goes through a case-by-case analysis, based on Lusztig's description of 
$\ph$-constructible characters and Lusztig $\ph$-families~\cite[\S{22}]{lusztig}. 

\bigskip

\noindent{\sc Remark 0 - } 
% Let $\gamma_\chi$ denote the coefficient of $\tb^{b_\chi}$ in 
% $F_\chi(\tb)$. Then it has been noticed by Lusztig~\cite[\S{2},~Page~325]{lusztig special} 
% that $\gamma_\chi=1$ whenever $\chi$ is special. 
%
As the only irreducible Coxeter systems affording possibly unequal parameters are of type 
$I_2(2m)$, $F_4$ or $B_n$, and as $r_{\! \chi}=1$ for any character $\chi$ in these groups, 
the statement ``$r_{\! \chi}=1$'' in Theorem~L(a) and (b) follows immediately 
from Lusztig's Theorem. Therefore, we will prove only the statements about the minimality 
of the $b$-invariant and the scalar product.\finl

\bigskip

\noindent{\bf Acknowledgements.} We wish to thank N. Jacon for pointing out a mistake in a preliminary 
version of this work.

\bigskip

\section{Proof of Theorem~L}

\medskip

\subsection{Reduction} 
It is easily seen that the proof of Theorem~L may be reduced to the case where 
$(W,S)$ is irreducible. If $W$ is of type $A_n$, $D_n$, $E_6$, $E_7$, $E_8$, $H_3$ or $H_4$, 
then $\ph$ is necessarily constant and Theorem~L follows immediately from 
Lusztig's Theorem. If $W$ is dihedral, then Theorem~L is easily checked 
using~\cite[\S{17}~and~Lemma~22.2]{lusztig}. 
If $W$ is of type $F_4$, then Theorem~L follows from inspection of~\cite[Table~2]{geck f4}. 
Therefore, this shows that we may, and we will, 
assume that $W$ is of type $B_n$, with $n \ge 2$. Write $S=\{t,s_1,s_2,\dots,s_{n-1}\}$ 
in such a way that the Dynkin diagram of $(W,S)$ is 
\begin{center}
\begin{picture}(250,30)
\put(-91,10){$(\#)$}
\put( 40, 10){\circle{10}}
\put( 44,  7){\line(1,0){33}}
\put( 44, 13){\line(1,0){33}}
\put( 81, 10){\circle{10}}
\put( 86, 10){\line(1,0){29}}
\put(120, 10){\circle{10}}
\put(125, 10){\line(1,0){20}}
\put(155,  7){$\cdot$}
\put(165,  7){$\cdot$}
\put(175,  7){$\cdot$}
\put(185, 10){\line(1,0){20}}
\put(210, 10){\circle{10}}
\put( 38, 20){$t$}
\put( 76, 20){$s_1$}
\put(116, 20){$s_2$}
\put(204, 20){$s_{n{-}1}$}
\end{picture}
\end{center}
Write $b=\ph(t)$ and $a=\ph(s_1)=\ph(s_2)=\cdots = \ph(s_{n-1})$. If $b \not\in a\NM^*$, 
then $\Cons_\ph^\kl(W)=\Irr(W)$ (see~\cite[Proposition~22.25]{lusztig}) 
and Theorem~L becomes obvious. So we may assume 
that $b =ra$ with $r \in \NM^*$, and since the notions are unchanged by multiplying $\ph$ 
by a positive real number, we may also assume that $a=1$. Therefore:

\bigskip

\begin{quotation}
{\bf Hypothesis and notation.} 
{\it From now on, and until the end of this section, we  assume that the Coxeter system 
$(W,S)$ is of type $B_n$, with $n \ge 2$, that $S=\{t,s_1,s_2,\dots,s_{n-1}\}$ is such that 
the Dynkin diagram of $(W,S)$ is given by $(\#)$ and that $\ph(t)=r\ph(s_1)=r\ph(s_2)=\cdots = r\ph(s_{n-1})=r$ 
with $r \in \NM^*$.
%  and that $\ph$ is not constant (i.e. that 
% $r \ge 2$).
}
\end{quotation}

\bigskip

We will now review the combinatorics introduced by Lusztig (symbols, admissible involutions,...) 
in order to compute families and constructible characters in type $B_n$ (see~\cite[\S{22}]{lusztig} 
for further details).

\bigskip

% 
% \bigskip
% We will now use the combinatorics of symbols, as Lusztig did for determining 
% the families and the constructible characters.

\def\symbols{{\mathbf{Sym}}}%\!\!m}

\subsection{Admissible involutions} 
Let $l \ge 0$ and let $Z$ be a totally ordered set of size $2l+r$. 
We will define by induction on $l$ what is an $r$-admissible involution of $Z$. Let 
$\iota : Z \to Z$ be an involution. Then $\iota$ is said {\it $r$-admissible} 
if it has $r$ fixed points and, if $l \ge 1$, there exist two consecutive elements 
$b$ and $c$ of $Z$ such that $\iota(b)=c$ and the restriction of $\iota$ to 
$Z \setminus \{b,c\}$ is $r$-admissible.

Note that, if $\iota$ is an $r$-admissible involution and if $\iota(b)=c > b$ and $\iota(z)=z$, 
then $z < b$ or $z > c$ (this is easily proved by induction on $|Z|$).

\bigskip

\subsection{Symbols} 
We will denote by $\symbols_k(r)$ the set of {\it symbols} $\L=\DS{\b \choose \g}$ where 
$\b=(\b_1 < \b_2 < \cdots < \b_{k+r})$ and $\g=(\g_1 < \g_2  < \cdots < \g_k)$ 
are increasing sequences of {\it non-zero} natural numbers. We set
$$|\L|=\sum_{i=1}^{k+r} (\b_i-i) + \sum_{j=1}^k (\g_j-j)$$
$$\bb(\L)=\sum_{i=1}^{k+r} (2k+2r-2i)(\b_i-i) + \sum_{j=1}^k (2k+1-2j)(\g_j-j).\leqno{\text{and}}$$
The number $\bb(\L)$ will be called the {\it $\bb$-invariant} of $\L$. For simplifying our arguments, 
we will define 
$$\nabla_{k,r}=\sum_{i=1}^{k+r} (2k+2r-2i)i + \sum_{j=1}^k (2k+1-2j)j$$
so that 
$$\bb(\L)=\sum_{i=1}^{k+r} (2k+2r-2i)\b_i + \sum_{j=1}^k (2k+1-2j)\g_j - \nabla_{k,r}.$$
% $$\bb'(\L)=\sum_{i=1}^{k+r} (2k+2r-2i)\b_i + \sum_{j=1}^k (2k+1-2j)\g_j=\bb(\L)+\nabla_{k,r}.
% \leqno{\text{and}}$$
By abuse of notation, we 
denote by $\b \cap \g$ the set $\{\b_1,\b_2,\dots,\b_{k+r}\} \cap \{\g_1,\g_2,\dots,\g_k\}$ 
and by $\b \cup \g$ the set $\{\b_1,\b_2,\dots,\b_{k+r}\} \cup \{\g_1,\g_2,\dots,\g_k\}$.
We also set $\b \dotplus \g = (\b \cup \g) \setminus (\b \cap \g)$.

We now define
$$\zb'(\L)=(\b_1,\b_2,\dots,\b_r,\g_1,\b_{r+1},\g_2,\b_{r+2},\dots,\g_k,\b_{r+k})$$
and we will write 
$$\zb'(\L)=(z_1'(\L),z_2'(\L),\cdots,z_{2k+r}'(\L)),$$
so that
$$\begin{array}{rcl}
\bb(\L)&=&\DS{\sum_{i=1}^r (2k+2r-2i)z_i'(\L) + \sum_{i=r+1}^{2k+r} (2k+r-i) z_i'(\L)} - \nabla_{k,r}\\
&=& \DS{\sum_{i=1}^r (r-i)z_i'(\L) + \sum_{i=1}^{2k+r} (2k+r-i) z_i'(\L)}  - \nabla_{k,r}\\
&=& \DS{\sum_{i=1}^{r-1} \Bigl(\sum_{j=1}^i z_j'(\L)\Bigr)  
+ \sum_{i=1}^{2k+r-1} \Bigl(\sum_{j=1}^i z_j'(\L)\Bigr)} - \nabla_{k,r}.\\
\end{array}\leqno{(\clubsuit)}$$ 

\bigskip

\subsection{Families of symbols} 
We denote by $\zb(\L)$ the sequence $z_1 \le z_2 \le \cdots \le z_{2k+r}$ obtained after 
rewriting the sequence $(\b_1,\b_2,\dots,\b_{k+r},\g_1,\g_2,\dots,\g_k)$ in 
non-decreasing order. 

\bigskip

% 
% \noindent{\bf Definition 1.} 
% {\it The symbol $\L$ is said to be {\bfit special} if $\zb(\L) = \zb'(\L)$.}
% 
% \bigskip

\noindent{\sc Remark 1 - } Note that the sequence $\zb'(\L)$ determines the symbol $\L$, contrarily to 
the sequence $\zb(\L)$. However, $\zb(\L)$ determines completely $|\L|$ thanks to the formula 
$|\L|=\sum_{z \in \zb(\L)} z - r(r+1)/2-(k+r)(k+r+1)/2$.~$\square$

\bigskip
% 
% \subsection{A comparison  of $\bb$-invariants of symbols} 
% Let $\L \in \symbols_k(r)$ and let $1 \le i < j \le 2k+r$. 
% We denote by $\L\{i,j\}$ the symbol

%\subsection{Families, constructible characters} 
We say that two symbols $\L=\DS{\b \choose \g}$ and $\L'=\DS{\b' \choose \g'}$ 
%where $\b'=(\b_1',\b_2',\dots,\b_{k+r}')$ and $\g'=(\g_1',\g_2',\dots,\g_{k}')$ 
in $\symbols_k(r)$ {\it are in the same family} if $\zb(\L)=\zb(\L')$. 
Note that this is equivalent to say that $\b\cap \g = \b' \cap \g'$ and $\b \cup \g = \b' \cup \g'$. 
If $\FC$ is the family of $\L$, we set $X_\FC=\b \cap \g$ and $Z_\FC=\b \dot{+} \g$: note that $X_\FC$ and 
$Z_\FC$ depend only on $\FC$ (and not on the particular choice of $\L \in \FC$). 

If $\iota$ is an $r$-admissible involution 
of $Z_\FC$, we denote by $\FC_\iota$ the set of symbols $\L=\DS{\b \choose \g}$ in $\FC$ such 
that $|\b \cap \o|=1$ for all $\iota$-orbits $\o$. 

\bigskip

\subsection{Lusztig families, constructible characters} 
Let $\L \in \symbols_k(r)$ be such that $|\L|=n$. Let $\Bip(n)$ be the set 
of bipartitions of $n$. We set
$$\lamb_1(\L)=(\b_{k+r}-(k+r) \ge \cdots \ge \b_2-2 \ge \b_1-1),$$
$$\lamb_2(\L)=(\g_k-k \ge \cdots \ge \g_2-2 \ge \g_1-1)$$
$$\lamb(\L)=(\lamb_1(\L),\lamb_2(\L)).\leqno{\text{and}}$$
Then $\lamb(\L)$ is a bipartition of $n$. We denote by $\chi_\L$ the irreducible 
character of $W$ denoted by $\chi_{\lamb(\L)}$ in~\cite[\S{22}]{lusztig} 
or in~\cite[\S{5.5.3}]{geck pfeiffer}. Then~\cite[\S{5.5.3}]{geck pfeiffer}
$$b_{\chi_\L}=\bb(\L).\leqno{(\diamondsuit)}$$
With these notations, Lusztig described the $\ph$-constructible characters 
in~\cite[Proposition~22.24]{lusztig}, from which the description of Lusztig $\ph$-families 
follows by using~\cite[Lemma~22.22]{lusztig}:

\bigskip

\noindent{\bf Theorem 2 (Lusztig).}
{\it Let $\FC_\kl$ be a Lusztig $\ph$-family and let $\g \in \Cons_\ph^\kl(\FC_\kl)$. 
If we choose $k$ sufficiently large, then:
\begin{itemize}
\itemth{a} There exists a family $\FC$ of symbols in $\symbols_k(r)$ such that
$$\FC_\kl=\{\chi_\L~|~\L \in \FC\}.$$

\itemth{b} There exists an $r$-admissible involution $\iota$ of $Z_\FC$ such that 
$$\g=\sum_{\L \in \FC_\iota} \chi_\L.$$
\end{itemize}}

% \bigskip
% 
% \noindent{\it Proof.} 
% (b) is proved in~\cite[Proposition~22.24]{lusztig} while (a) 
% follows from (b) and~\cite[Lemma~22.22]{lusztig}.\fin

\bigskip

% If $z \in Z_\FC$ and $\L \in \FC$, we denote by $\Pos(z,\L)$ the unique natural number $i$ 
% such that $z=z_i'(\L)$: it is unique because $z \not\in X_\FC$.
% 
% \bigskip

If $\L=\DS{\b \choose \g}$, we set $\L^\#=\DS{\b \setminus (\b \cap \g) \choose \g \setminus (\b \cap \g)}$. 

\bigskip

\noindent{\bf Definition 3.} 
{\it The symbol $\L$ is said {\bfit special} if $\zb(\L^\#)=\zb'(\L^\#)$. 
% If $\L \in \FC_\iota$, then $\L$ is called 
% {\bfit $\iota$-semispecial} if the following two properties are satisfied:
% \begin{quotation}
% \begin{itemize}
% \itemth{S1} If $\o$ is a $\iota$-orbit, then $\Pos(\min(\o),\L) \le \Pos(\max(\o),\L)$. 
% 
% \itemth{S2} If $\o$ and $\o'$ are two $\iota$-orbits such that $\min(\o) < \min(\o')$, 
% then $\Pos(\min(\o),\L) < \Pos(\min(\o'),\L)$.
% \end{itemize}
% \end{quotation}
}

\bigskip

\noindent{\sc Remark 4 - } 
According to Remark 1, there is a unique special symbol in each family. It will be denoted by 
$\L_\FC$. Finally, note that, if $\L$, $\L'$ belong to the same family, then $|\L|=|\L'|$.~$\square$ 

\bigskip

Now, Theorem~L follows from Theorem~2, Formula~$(\diamondsuit)$ and the following next Theorem:

\bigskip

\noindent{\bf Theorem 5.} {\it Let $\FC$ be a family of symbols in $\symbols_k(r)$, 
let $\iota$ be an $r$-admissible involution of $Z_\FC$ and let $\L \in \FC$. Then:
\begin{itemize}
\itemth{a} $\bb(\L) \ge \bb(\L_\FC)$ with equality if and only if $\L=\L_\FC$.

\itemth{b} There is a unique symbol $\L_{\FC,\iota}$ in $\FC_\iota$ such that, 
if $\L \in \FC_\iota$, then $\bb(\L) \ge \bb(\L_{\FC,\iota})$, with equality 
if and only if $\L=\L_{\FC,\iota}$. 
\end{itemize}}

\bigskip

The rest of this section is devoted to the proof of Theorem~5.

\bigskip

\subsection{First reduction} 
First, assume that $X_\FC \neq \vide$. Let $b \in X_\FC$ and let $\L = \DS{\b \choose \g}\in \FC$.
Then $b \in \b \cap \g = X_\FC$ and we denote by $\b[b]$ the sequence obtained by removing $b$ to $\b$. 
Similarly, let $\L[b] = \DS{\b[b] \choose \g[b]}$. 

Then $\L[b] \in \symbols_{k-1}(r)$ and 

$$\bb(\L)=\bb(\L[b])+ \nabla_{k,r}-\nabla_{k-1,r} + 
b \Bigl(4k+2r+1-\sum_{\stackrel{\SS{z \in \zb(\L)}}{z \le b}} 2\Bigr) 
+ 2 \sum_{\stackrel{\SS{z \in \zb(\L)}}{z < b}} z.\leqno{(\heartsuit)}$$

\medskip

% \begin{quotation}
% {\small 
\noindent{\it Proof of $(\heartsuit)$.} 
Let $i_0$ and $j_0$ be such that $\b_{i_0} = b$ and $\g_{j_0}=b$. Then
%\eqna
$$\bb(\L)-\bb(\L[b]) = \nabla_{k,r}-\nabla_{k-1,r} \DS{+ 
(2k+2r-2i_0)b + \sum_{i=1}^{i_0-1} 2 \b_i + (2k+1-2j_0) b + \sum_{j=1}^{j_0-1} 2 \g_j.}
$$
But the numbers $\b_1$, $\b_2$,\dots, $\b_{i_0}$, $\g_1$, $\g_2$,\dots, $\g_{j_0}$ 
are exactly the elements of the sequence $\zb(\L)$ which are $\le b$. 
So 
$$i_0+j_0 = \sum_{\stackrel{\SS{z \in \zb(\L)}}{z \le b}} 1$$
and
$$\sum_{i=1}^{i_0-1} \b_i + \sum_{j=1}^{j_0-1} \g_j = \sum_{\stackrel{\SS{z \in \zb(\L)}}{z < b}} z.$$
This shows $(\heartsuit)$.\fin
% }
% \end{quotation}

\bigskip

 Now, the family of $\L[b]$ depends only on the family of $\L$ (and not on $\L$ itself): 
indeed, $\zb(\L[b])$ is obtained from $\zb(\L)$ by removing the two entries equal to $b$. 
We will denote by $\FC[b]$ the family of $\L[b]$. Moreover, $Z_{\FC[b]}=Z_\FC$ and the map 
$\L \mapsto \L[b]$ induces a bijection between $\FC$ and $\FC[b]$, and also 
induces a bijection between $\FC_\iota$ and $\FC[b]_\iota$. 

On the other hand, the formula $(\heartsuit)$ 
shows that the difference between $\bb(\L)$ and $\bb(\L[b])$ depends only on $b$ 
and $\FC$, so proving Theorem~5 for the pair $(\FC,\iota)$ is equivalent to 
proving Theorem~5 for the pair $(\FC[b],\iota)$. By applying several 
times this principle if necessary, this means that 
we may, and we will, assume that
$$X_\FC = \vide.$$%\leqno{({\mathrm{Hyp}})}$$

\bigskip

\subsection{Proof of Theorem 5(a)} 
% Write $\zb'(\L)=(z_1',z_2',\cdots,z_{2k+r}')$ and 
% $\zb(\L)=\zb(\L_\FC)=\zb'(\L_\FC)=(z_1,z_2,\cdots,z_{2k+r})$. 
First, note that $\zb(\L)=\zb(\L_\FC)=\zb'(\L_\FC)$ (the last equality follows from the fact that 
$\L_\FC$ is special and $X_\FC=\vide$). 
As $\zb'(\L)$ is a permutation 
of the non-decreasing sequence $\zb'(\L_\FC)$, we have
$$\sum_{j=1}^i z_j'(\L) \ge \sum_{j=1}^i z_j'(\L_\FC)$$
for all $i \in \{1,2,\cdots,2k+r\}$. 
%, and
%$$\sum_{i=1}^{2k+r} z_i' = \sum_{i=1}^{2k+r} z_i$$. 
So, it follows from $(\clubsuit)$ that 
$$\bb(\L)-\bb(\L_\FC)=\sum_{i=1}^{r-1} \Bigl(\sum_{j=1}^i \bigl(z_j'(\L)-z_j'(\L_\FC)\bigr)\Bigr) + 
\sum_{i=1}^{2k+r-1} \Bigl(\sum_{j=1}^i \bigl(z_j'(\L)-z_j'(\L_\FC)\bigr)\Bigr).$$
So $\bb(\L) \ge \bb(\L_\FC)$ with equality only whenever 
$\sum_{j=1}^i z_j'(\L)=\sum_{j=1}^i z_j'(\L_\FC)$ for all $i \in \{1,2,\dots,2k+r\}$. 
The proof of Theorem 5(a) is complete.

% 
% \subsection{Operations on symbols and $\bb$-invariant} 
% We fix in this subsection a symbol $\L=\DS{\b \choose \g} \in \symbols_k(r)$. 
% 
% \medskip
% 
% 
% 
% 
\bigskip

\subsection{Proof of Theorem 5(b)} 
We denote by $f_r < \cdots <f_1$ the elements of $Z_\FC$ which are fixed by $\iota$. 
We also set $f_{r+1} = 0$ and $f_0 = \infty$. 
As $\iota$ is $r$-admissible, the set $Z_\FC^{(d)}=\{z \in Z_\FC~|~f_{d+1} < z < f_d\}$ is $\iota$-stable 
and contains no $\iota$-fixed point (for $d \in \{0,1,\dots,r\}$). 
Let $k_d=|Z_\FC^{(d)}|/2$ and let $\iota_d$ be the restriction 
of $\iota$ to $Z_\FC^{(d)}$. Then $\iota_d$ is a $0$-admissible involution of $Z_\FC^{(d)}$. 

If $\L=\DS{\b \choose \g} \in \FC_\iota$, we set 
$\b^{(d)} = \b \cap Z_\FC^{(d)}$, $\g^{(d)} = \g \cap Z_\FC^{(d)}$ and $\L^{(d)} = \DS{\b^{(d)} \choose \g^{(d)}}$. 
Then $\L^{(d)} \in \symbols_{k_d}(0)$ and, if $\FC^{(d)}$ denotes the family of $\L^{(d)}$, 
then $\L^{(d)} \in \FC^{(d)}_{\iota_d}$. 

Now, if $\L' = {\DS{\b' \choose \g'}} \in \symbols_{k'}(0)$, we set
$$\bb_d(\L') = \sum_{i=1}^{k'} (2k'+2d-2i) \b_i' + \sum_{j=1}^{k'} (2k'+1-2j) \g_j'.$$
The number $\bb_d(\L')$ is called the {\it $\bb_d$-invariant} of $\L'$. 
It then follows from the definition of $\bb$ and $\nabla_{k,r}$ that
$$
\bb(\L)=\DS{\sum_{d=0}^r \bb_d(\L^{(d)}) - \nabla_{k,r} 
+ \sum_{d=1}^r 2\bigl(k_0+k_1+\cdots + k_{d-1}\bigr) \Bigl(f_d + \sum_{z \in Z^{(d)}} z\Bigr).} 
\leqno{(\spadesuit)}$$
Since the map 
$$\fonctio{\FC_\iota}{\prod_{d=0}^r \FC_{\iota_d}^{(d)}}{\L}{(\L^{(0)},\L^{(1)},\dots,\L^{(d)})}$$
is bijective and since $\bb(\L) - \sum_{d=0}^r \bb_d(\L^{(d)})$ depends only on $(\FC,\iota)$ 
and not on $\L$ (as shown by the formula~$(\spadesuit)$), Theorem~5(b) will follow from the following 
lemma~:

\bigskip

\noindent{\bf Lemma 6.} 
{\it There exists a unique symbol in $\FC_{\iota_d}^{(d)}$ with minimal $\bb_d$-invariant.}

\bigskip

The proof of Lemma~6 will be given in the next section.

\bigskip

\section{Minimal $\bb_d$-invariant} 

\medskip

For simplifying notation, we set $Z=Z_\FC^{(d)}$, $l=k_d$, $\GC=\FC^{(d)}$ and $\jmath=\iota_d$. 
Let us write $Z=\{z_1,z_2,\dots,z_{2l}\}$ with $z_1 < z_2 < \cdots < z_{2l}$. Recall from the 
previous secion that $\jmath$ 
is a $0$-admissible involution of $Z$. 

\bigskip

\subsection{Construction} 
We will define by induction on $l \ge 0$ a symbol $\L_\jmath^{(d)}(Z) \in \GC_\jmath$. 
If $l=0$, then $\L_\jmath^{(d)}(Z)$ is obviously empty. So assume now that, for any 
set of non-zero integers $Z'$ of order $2(l-1)$, for any $0$-admissible involution $\jmath'$ 
of $Z'$ and any $d' \ge 0$, we have defined a symbol $\L_{\jmath'}^{(d')}(Z')$. 
Then $\L^{(d)}_\jmath(Z)=\DS{\b_\jmath^{(d)}(Z) \choose \g_\jmath^{(d)}(Z)}$ 
is defined as follows: let $Z' = Z\setminus\{z_1,\iota(z_1)\}$, 
$\jmath'$ the restriction of $\jmath$ to $Z'$ and let 
$$d'=
\begin{cases}
d-1 & \text{if $d \ge 1$,}\\
1 & \text{if $d=0$.}
\end{cases}
$$
Then $|Z'|=2(l-1)$ and $\jmath'$ is $0$-admissible. So $\L_{\jmath'}^{(d')}(Z')=
\DS{\b_{\jmath'}^{(d')}(Z') \choose \g_{\jmath'}^{(d')}(Z')}$ is 
well-defined by the induction hypothesis. We then set 
$$\b_\jmath^{(d)}(Z)=
\begin{cases}
\b_{\jmath'}^{(d')}(Z') \cup \{z_1\} & \text{if $d \ge 1$,}\\
\b_{\jmath'}^{(d')}(Z') \cup \{\jmath(z_1)\} & \text{if $d=0$,}
\end{cases}
$$
$$\g_\jmath^{(d)}(Z)=
\begin{cases}
\g_{\jmath'}^{(d')}(Z') \cup \{\jmath(z_1)\} & \text{if $d \ge 1$,}\\
\g_{\jmath'}^{(d')}(Z') \cup \{z_1\} & \text{if $d=0$.}
\end{cases}\leqno{\text{and}}
$$
Then Lemma~6 is implied by the next lemma~:

\bigskip

\noindent{\bf Lemma ${\boldsymbol{6^+}}$.} 
{\it Let $\L \in \GC_\jmath$. Then $\bb_d(\L) \ge \bb_d(\L_\jmath^{(d)}(Z))$ with equality 
if and only if $\L=\L_\jmath^{(d)}(Z)$.}

\bigskip

The rest of this section is devoted to the proof of Lemma~$6^+$. 
We will first prove Lemma~$6^+$ whenever $d \in \{0,1\}$ using Lusztig's Theorem. 
We will then turn to the general case, which will be handled by induction on $l=|Z|/2$. 
We fix $\L= {\DS{\b \choose \g}} \in \GC_\iota$.

\def\dotcup{\hskip1mm\dot{\cup}\hskip1mm}
\def\oppose{{\mathrm{op}}}

\bigskip

\subsection{Proof of Lemma ${\boldsymbol{6^+}}$ whenever ${\boldsymbol{d=1}}$} 
Let $z$ be a natural number strictly bigger than all the elements of $Z$. 
Let $\Lamt=\DS{\b \cup \{z\} \choose \g} \in \symbols_k(1)$. Then $\bb_1(\L)=\bb(\Lamt) + C$, where 
$C$ depends only on $Z$. Let $\Lamt_0=\DS{z_1 ,z_3,\dots,z_{2l-1},z \choose z_2,\dots,z_{2l}}$. 
Since $\jmath$ is $0$-admissible, it is easily seen that, if $\jmath(z_i)=z_j$, then $j-i$ is odd. 
So $\Lamt_0 \in \GC_\jmath$. But, by~\cite[\S{5}]{lusztig special}, 
$\bb(\Lamt) \ge \bb(\Lamt_0)$ with equality if and only if $\Lamt=\Lamt_0$. 
So it is sufficient to notice that $\widetilde{\L_\jmath^{(1)}(Z)}=\Lamt_0$, 
which is easily checked.

\bigskip

\subsection{Proof of Lemma ${\boldsymbol{6^+}}$ whenever ${\boldsymbol{d=0}}$} 
Assume in this subsection, and only in this subsection, that $d=0$ or $1$. We denote by 
$\L^\oppose=\DS{\g \choose \b} \in \GC_\jmath$. It is readily seen from the construction that 
$\L_\jmath^{(0)}(Z)^\oppose=\L_\jmath^{(1)}(Z)$ and that 
$$\bb_1(\L)=\bb_0(\L^\oppose) + \sum_{z \in Z} z.$$
So Lemma~$6^+$ for $d=0$ follows from Lemma~$6^+$ for $d=1$.

\bigskip

\subsection{Proof of Lemma ${\boldsymbol{6^+}}$ whenever ${\boldsymbol{d \ge 2}}$} 
Assume now, and until the end of this section, that $d \ge 2$. 
We will prove Lemma~$6^+$ by induction on $l=|Z|/2$. The result is obvious if $l=0$, as well as 
if $l = 1$. So we assume that $l \ge 2$ and that Lemma $6^+$ holds for $l' \le l-1$. Write $\jmath(z_1)=z_{2m}$, 
where $m \le l$ (note that $\jmath(z_1) \not\in \{z_1,z_3,z_5,\dots,z_{2l-1}\}$ since $\jmath$ is 
$0$-admissible). 

\bigskip

Assume first that $ m < l$. Then $Z$ can we written as the union 
$Z=Z^+ \dotcup Z^-$, where $Z^+=\{z_1,z_2,\dots,z_{2m}\}$ and $Z^-=\{z_{2m+1},z_{2m+2},\dots,z_{2l}\}$ 
are $\jmath$-stable (since $\jmath$ is $0$-admissible). If $\e \in \{+,-\}$, let $\jmath^\e$ 
denote the restriction of $\jmath$ to $Z^\e$, let $\b^\e=\b \cap Z^\e$, $\g^\e = \g \cap Z^\e$ and  
$\L^\e=\DS{\b^\e \choose \g^\e}$, and let $\GC^\e$ denote the family of $\L^\e$. Then it is easily seen that 
$\L^\e \in \GC^\e_{\jmath^\e}$, that $\bb_d(\L)-\bigl(\bb_d(\L^+)+\bb_d(\L^-)\bigr)$ depends only on $(\GC,\jmath)$ 
and that $\L_\jmath^{(d)}(Z)^\e=\L_{\jmath^\e}^{(d)}(Z^\e)$. 
By the induction hypothesis, $\bb_d(\L^\e) \ge \bb_d(\L_{\jmath^\e}^{(d)}(Z^\e))$ with 
equality if and only if $\L^\e=\L_{\jmath^\e}^{(d)}(Z^\e)$. So the result follows in this case. 
This means that we may, and we will, work under the following hypothesis:

\bigskip

\boitegrise{\noindent{\bf Hypothesis.} 
{\it From now on, and until the end of this section, we assume that $\jmath(z_1)=z_{2l}$.}}{0.75\textwidth}

\bigskip

As in the construction of $\L_\jmath^{(d)}(Z)$, let $Z' = Z\setminus\{z_1,z_{2l}\}
=\{z_2,z_3,\dots,z_{2l-1}\}$, let $\jmath'$ denote 
the restriction of $\jmath$ to $Z'$ and let 
$$d'=
\begin{cases}
d-1 & \text{if $d \ge 1$,}\\
1 & \text{if $d=0$.}
\end{cases}
$$
Then $|Z'|=2(l-1)$ and $\jmath'$ is $0$-admissible. Let $\L'=\DS{\b' \choose \g'}$ where 
$\b'=\b \setminus \{z_1,z_{2l}\}$ and $\g'=\g\setminus\{z_1,z_{2l}\}$. Since $d \ge 2$, 
we have $z_1 \in \b_\jmath^{(d)}(Z)$ and $z_{2l} \in \g_\jmath^{(d)}(Z)$. 
This implies that
$$\bb_d(\L_\jmath^{(d)}(Z)) = \bb_{d-1}(\L_{\jmath'}^{(d-1)}(Z')) + 
z_{2l} + 2(l+d) z_1 + 2 \sum_{z \in Z'} z.\leqno{(\bigstar)}$$
If $z_1 \in \b$, then $\L=\L_\jmath^{(d)}(Z)$ if and only if $\L'=\L_{\jmath'}^{(d')}(Z')$ 
and again 
$$\bb_d(\L) = \bb_{d-1}(\L') + 
z_{2l} + 2(l+d) z_1 + 2 \sum_{z \in Z'} z.$$
So the result follows from $(\bigstar)$ and from the induction hypothesis. 

This means that we may, and we will, assume that $z_1 \in \g$. In this case,
$$\bb_d(\L) = \bb_{d+1}(\L') + 2d z_{2l}+ (2l+1) z_1.$$
Then it follows from $(\bigstar)$ that 
$$\bb_d(\L) - \bb_d(\L_\jmath^{(d)}(Z)) = \bb_{d+1}(\L') - \bb_{d-1}(\L_{\jmath'}^{(d-1)}(Z')) 
+ (2d-1) (z_{2l}-z_1) - 2 \sum_{z \in Z'} z.$$
So, by the induction hypothesis, 
$$\bb_d(\L) - \bb_d(\L_\jmath^{(d)}(Z)) \ge \bb_{d+1}(\L_{\jmath'}^{(d+1)}(Z')) - \bb_{d-1}(\L_{\jmath'}^{(d-1)}(Z')) 
+ (2d-1) (z_{2l}-z_1) - 2 \sum_{z \in Z'} z.$$
Since $z_{2l}-z_1 > z_{2l-1} - z_2$, it is sufficient to show that 
$$\bb_{d+1}(\L_{\jmath'}^{(d+1)}(Z')) - \bb_{d-1}(\L_{\jmath'}^{(d-1)}(Z')) 
\ge - (2d-1) (z_{2l-1}-z_2) + 2 \sum_{z \in Z'} z.\leqno{(?)}$$
This will be proved by induction on the size of $Z'$. 
First, if $\jmath(z_2) < z_{2l-1}$, then we can separate $Z'$ into two $\jmath'$-stable 
subsets and a similar argument as before allows to conclude thanks to the induction hypothesis. 

So we assume that $\jmath'(z_2)=z_{2l-1}$. Let $Z''=Z'\setminus\{z_2,z_{2l-1}\}$ 
and let $\jmath''$ denote the restriction of $\jmath'$ to $Z''$. Since $z_2\in\b_{\jmath'}^{(d+1)}(Z')$, 
we can apply $(\bigstar)$ one step further to get
\eqna
\bb_{d+1}(\L_{\jmath'}^{(d+1)}(Z')) - \bb_{d-1}(\L_{\jmath'}^{(d-1)}(Z')) 
&\!\!=\!\!& \DS{\bb_d(\L_{\jmath''}^{(d)}(Z'') + z_{2l-1} + 2(l+d)z_2 + 2 \sum_{z \in Z''} z }\\
&& \DS{- \bigl( \bb_{d-2}(\L_{\jmath''}^{(d-2)}(Z'')) + z_{2l-1} + 2(l+d-2) z_2 + 2 \sum_{z \in Z''} z\bigr)}\\
&\!\!=\!\!&\bb_d(\L_{\jmath''}^{(d)}(Z''))-\bb_{d-2}(\L_{\jmath''}^{(d-2)}(Z'')) + 4 z_2. \\
\endeqna
So, by the induction hypothesis,
\eqna
\bb_{d+1}(\L_{\jmath'}^{(d+1)}(Z')) - \bb_{d-1}(\L_{\jmath'}^{(d-1)}(Z')) 
&\ge&\DS{-(2d-3) (z_{2l-2}-z_3) + 2 \sum_{z \in Z''} z + 4 z_2}\\
&\ge& \DS{-(2d-3) (z_{2l-1}-z_2) + 2 \sum_{z \in Z'} z + 2 z_2 - 2 z_{2l-1}}\\
&=&\DS{-(2d-1) (z_{2l-1}-z_2) + 2 \sum_{z \in Z'} z },\\
\endeqna
as desired. This shows $(?)$ and completes the proof of Lemma~$6^+$.

\bigskip

\section{Complex reflection groups}

\medskip

If $\WC$ is a complex reflection group, then R. Rouquier and the author have also defined 
{\it Calogero-Moser cellular characters} and {\it Calogero-Moser families} (see~\cite{note} 
or~\cite{cmcells}). If $\WC$ is of type $G(l,1,n)$ (in Shephard-Todd classification), 
then Leclerc and Miyachi~\cite[\S{6.3}]{leclerc} proposed, in link with canonical bases 
of $U_v(\sG\lG_\infty)$-modules, a family of characters that could be a good analogue 
of constructible characters: let us call them the {\it Leclerc-Miyachi constructible 
characters} of $G(l,1,n)$. If $l=2$, then they coincide with 
constructible characters~\cite[Theorem~10]{leclerc}. 

Of course, it would be interesting to know if Calogero-Moser cellular 
characters coincide with the Leclerc-Miyachi ones: this seems rather 
complicated but it should be at least possible to check if 
the Leclerc-Miyachi constructible characters satisfy 
the analogous properties with respect to the $b$-invariant.


\bigskip

\begin{thebibliography}{AAAA}
\bibitem[Be1]{bellamy these} {\sc G. Bellamy}, 
Generalized Calogero-Moser spaces and rational Cherednik algebras, 
PhD thesis, University of Edinburgh (2010).

\bibitem[Be2]{bellamy} {\sc G. Bellamy}, 
The Calogero-Moser partition for $G(m,d,n)$, {\it Nagoya Math. J.} {\bf 207} (2012), 47-77.


\bibitem[BoRo1]{note} {\sc C. Bonnaf\'e \& R. Rouquier}, Calogero-Moser versus Kazhdan-Lusztig cells, 
{\it Pacific J. Math.} {\bf 261} (2013), 45-51.

\bibitem[BoRo2]{cmcells} 
{\sc C. Bonnaf\'e \& R. Rouquier}, Cellules de Calogero-Moser, 
preprint (2013), {\tt arXiv:1302.2720}. 

\bibitem[Ge]{geck f4} {\sc M. Geck}, 
Computing Kazhdan-Lusztig cells for unequal parameters
{\it J. Algebra} {\bf 281} (2004) 342-365.
% 
% \bibitem[Ge2]{geck pycox} {\sc M. Geck}, 
% PyCox: Computing with (finite) Coxeter groups and Iwahori-Hecke algebras,
% To appear in {\it London Math. Soc. J. of Comput. and Math.}, preprint (2012), 
% {\tt arXiv:1201.5566}.

\bibitem[GePf]{geck pfeiffer} 
{\sc M. Geck \& G. Pfeiffer}, 
{\it Characters of finite Coxeter groups and Iwahori-Hecke algebras}, 
London Mathematical Society Monographs, New Series {\bf 21}, The Clarendon Press, 
Oxford University Press, New York, 2000, xvi+446 pp.

\bibitem[GoMa]{gordon martino} {\sc I. G. Gordon \& M. Martino}, 
Calogero-Moser space, restricted rational Cherednik algebras and two-sided cells, 
{\it Math. Res. Lett.} {\bf 16} (2009), 255-262. 

\bibitem[LeMi]{leclerc} {\sc B. Leclerc \& H. Miyachi}, 
Constructible characters and canonical bases, 
{\it J. Algebra} {\bf 277} (2004), 298-317.
 
\bibitem[Lu1]{lusztig special} {\sc G. Lusztig}, A class of irreducible representations 
of a Weyl group, {\it Indag. Math.} {\bf 41} (1979), 323-335.

\bibitem[Lu2]{lusztig orange} {\sc G. Lusztig}, 
{\it Characters of reductive groups over a finite field}, 
Annals of Mathematical Studies, {\bf 107}. 
Princeton University Press, Princeton, NJ, 1984, 
xxi+384 pp.

\bibitem[Lu3]{lusztig} 
{\sc G.~Lusztig}, 
\emph{Hecke algebras with unequal parameters}, CRM Monograph Series 
{\bf 18}, American Mathematical Society, Providence, RI (2003), 136 pp.

\end{thebibliography}
\end{document}